\documentclass[11pt,reqno]{amsart}
\usepackage{fullpage,amsfonts,amsmath,amscd,amssymb}

\newtheorem*{lemma}{Lemma}

\newtheorem*{theorem}{Theorem}
\newtheorem*{corollary}{Corollary}

\newtheorem*{prop}{Proposition}

{
}

\theoremstyle{definition}

\newtheorem*{question}{Question}

\theoremstyle{remark}

\newtheorem*{remark}{Remark}

\newcommand{\hc}{H_{\mathbf{c}}}
\newcommand{\xc}{X_{\mathbf{c}}}
\newcommand{\zc}{Z_{\mathbf{c}}}
\newcommand{\yc}{Y_{\mathbf{c}}}

\newcommand{\psiol}{{\overline \psi}}
\newcommand{\Rol}{{\overline R}}

\def\CC{{\mathbb C}}
\def\ZZ{{\mathbb Z}}
\def\PP{{\mathbb P}}
\def\RR{{\mathbb R}}

\def\fm{{\mathfrak{m}}}
\def\fp{{\mathfrak{p}}}
\def\fq{{\mathfrak{ q}}}

\DeclareMathOperator{\Spec}{Spec}

\DeclareMathOperator{\coh}{coh} \DeclareMathOperator{\spec}{Spec}

\DeclareMathOperator{\SL}{SL}
\DeclareMathOperator{\mat}{Mat} 
\DeclareMathOperator{\edo}{End} 
\DeclareMathOperator{\End}{End}
\DeclareMathOperator{\Hom}{Hom}
\DeclareMathOperator{\Sets}{Sets}

\DeclareMathOperator{\Kdim}{Kdim}
\DeclareMathOperator{\gldim}{gldim}
\DeclareMathOperator{\rank}{rank}
\DeclareMathOperator{\hgt}{ht}
\DeclareMathOperator{\sym}{Sym} 
\DeclareMathOperator{\homo}{Hom} 
 \DeclareMathOperator{\id}{Id}
 
 \DeclareMathOperator{\md}{mod}
 \DeclareMathOperator{\hil}{Hilb}
\DeclareMathOperator{\gr}{gr} \DeclareMathOperator{\codim}{codim}

\begin{document}

\pagenumbering{arabic}

\title[symplectic reflection algebras]{Representations of symplectic reflection algebras and resolutions of deformations of symplectic quotient singularities}

\author{Iain Gordon and S. Paul Smith}
\address{Department of Mathematics, University of Glasgow, Glasgow G12 8QW}
\email{ig@maths.gla.ac.uk}
\address{ Department of Mathematics, Box 354350, University of Washington, Seattle, WA 98195}
\email{smith@math.washington.edu}

\thanks{The first author was partially supported by the Nuffield Foundation
grant NAL/00625/G and by the University of Washington's Milliman
Fund. The second author was partially supported by an NSF grant
DMS-0070560. Both authors are grateful to the Edinburgh
Mathematical Society and the Leverhulme Trust for support.}
\begin{abstract}
We give an equivalence of triangulated categories between the
derived category of finitely generated representations of
symplectic reflection algebras associated with wreath products
(with parameter $t=0$) and the derived category of coherent
sheaves on a crepant resolution of the spectrum of the centre of
these algebras.
\end{abstract}
\maketitle

\section{Introduction}
\label{intro}
\subsection{}
In this paper we take a step towards a geometric understanding of the representation theory of certain symplectic reflection algebras (with parameter $t=0$). A number of papers have shown that this representation theory is closely related to the singularities of the centre of these algebras and to resolutions of these singularities, \cite{EG}, \cite{BG}, \cite{GK}. Here we make such a relationship precise by proving that the category of finitely generated modules is derived equivalent to the category of coherent sheaves on an appropriate desingularisation.

 A long term goal in this project is to find character formulae for simple modules, generalising the work \cite{FG} and \cite{gor} in which Kostka polynomials appear. A simple consequence of the derived equivalence is a geometric interpretation for the number of simple modules of a symplectic reflection algebra with given central character. A closer analysis will undoubtedly reveal more.

\subsection{}
Let us summarise our results.
Let $\Gamma$ be a non-trivial finite subgroup of $\SL(2,\CC)$ and $n$
a positive integer.
 Let $\hc$ be the symplectic reflection algebra (with parameter $t=0$) for the wreath product $\Gamma_n = S_n\ltimes \Gamma^n$. (All undefined notation and definitions can be found later in the paper.) The spectrum of the centre of this algebra, $\xc = \spec \zc$, is a deformation of the symplectic quotient singularity $\CC^{2n}/\Gamma_n$. Our principal
result is the following.

\begin{theorem}
\label{main.thm}
There is a crepant resolution $\pi_{\mathbf{c}} : Y_\mathbf{c}
\longrightarrow \xc$ such that there is an equivalence of
triangulated categories
$$D^b(\md \hc) \longrightarrow D^b(\coh Y_{\mathbf{c}}),$$
between the bounded derived category of finitely generated $\hc$--modules and the bounded derived category of coherent sheaves on $Y_\mathbf{c}$. \end{theorem}

\subsection{}
In the special case $\mathbf{c}=\mathbf{0}$ the variety
$X_{\mathbf{c}}$ is the orbit space $\mathbb{C}^{2n}/\Gamma_n$, so
the above theorem includes the results of \cite{KV} on Kleinian
singularities ($n=1$) and the observation of \cite[Section
4.4]{wang1} (for general $n$). The proof we give here, however, is
by deformation from the $\mathbf{c}=\mathbf{0}$ case and so
depends on these results. To prove the equivalence, we use the
methods of \cite{BKR}, which were adapted to a mildly
non--commutative situation in \cite{vdb2}, together with results
from \cite{hai} and \cite{hai2}.

\subsection{}
For $x \in \xc$, let $\mathfrak{m}_x$ be the corresponding maximal ideal of $\zc$. The simple modules of the finite dimensional algebra  $\hc/\mathfrak{m}_x\hc$ are the simple $\hc$--modules  with central character $x$. The following corollary is a straightforward consequence of the above theorem.

\begin{corollary}
Let $\pi_{\mathbf{c}} :Y_{\mathbf{c}}\longrightarrow \xc$ be the crepant resolution in Theorem \ref{main.thm}. For all $x\in \xc$, there is an isomorphism of Grothendieck groups $K(\hc/\mathfrak{m}_x\hc) \cong K(\pi_{\mathbf{c}}^{-1}(x))$.
\end{corollary}

\subsection{}
In the case $n=1$ this recovers known results on the simple modules of deformed preprojective algebras, whilst for $\mathbf{c}=0$ the content is essentially the (generalised) McKay correspondence, \cite{kal}.

\subsection{} Symplectic reflection algebras also exist for finite Coxeter groups, $W$.
However, \cite[Theorem 1.1]{GK} shows that the only the orbits
spaces $\mathbb{C}^{2n}/W$ admitting crepant resolutions are for
$W$ of type $A$ and $B$, that is $W = S_n$ or $W = S_n\ltimes
(\mathbb{Z}/2\mathbb{Z})^n$. Thus there is no analogue of Theorem
\ref{main.thm} valid for all finite Coxeter groups; our result is
as general as we can expect.

\subsection{}
The paper is organised as follows. In Section \ref{sra} we recall
the definition and basic properties of symplectic reflection
algebras. A discussion of non--commutative crepant resolutions and
derived equivalences is given in Section \ref{nc}. Section
\ref{ss} presents some results on deformations of semi--small
morphisms and their relation to crepant resolutions and derived
equivalences. Finally, in Section \ref{main}, we prove our main
result and discuss the application to counting simple modules.

\section{Symplectic reflection algebras}
\label{sra}
\subsection{}
Let $\tilde{\omega}$ be the standard symplectic form on
$\mathbb{C}^2$, $\Gamma$ a finite subgroup of $SL(2,\mathbb{C})$
and $n$ a positive integer. The wreath product $\Gamma_n
\equiv S_n\ltimes \Gamma^n$ acts on $V\equiv (\mathbb{C}^2)^n$,
preserving the symplectic form $\omega \equiv \tilde{\omega}^n$.

\subsection{} Recall that $\gamma \in \Gamma_n$
acting on $V$ is called a \textit{symplectic
reflection} if $\dim (1-\gamma)(V) = 2$. The set of all symplectic
reflections is denoted by $\mathcal{S}$. Let $\mathbf{c}$ be a
$\CC$--valued function on $\mathcal{S}$, constant on conjugacy
classes ($\gamma \mapsto c_{\gamma}$). Given $\gamma\in
\mathcal{S}$ define the form $\omega_{\gamma}$ on $V$ to have radical
$\ker (1-\gamma)$ and to be the restriction of $\omega$ on the
$(1-\gamma)(V)$.
\subsection{} The \textit{symplectic reflection algebra} $\hc$ is
the $\CC$--algebra, defined as the quotient of the skew group ring $TV\ast \Gamma_n$
by the relations
$$ x\otimes y - y\otimes x = \sum_{\gamma\in \mathcal{S}}
c_{\gamma} \omega_{\gamma}(x,y)\gamma,$$ for all $x,y\in V$.

\begin{remark} Usually, symplectic reflection algebras depend on a
further parameter $t\in \CC$, \cite[Section 1]{EG}. The definition above is
the case $t=0$.
\end{remark}

\subsection{} \label{ff} There is an increasing $\mathbb{N}$--filtration on
$\hc$, obtained by setting $F^0\hc = \CC\Gamma_n$, $F^1 = V\otimes
\CC\Gamma_n + \CC\Gamma_n$, and $F^i = (F^1)^i$. By the PBW
theorem, \cite[Theorem 1.3]{EG}, $\gr \hc \cong \CC[V]\ast \Gamma_n$ (where we have used
$\omega$ to identify the $\Gamma_n$--spaces $V$ and $V^*$). In
particular, a non--zero $\mathbf{c}$ yields a flat family of
symplectic reflection algebras $H_{u\mathbf{c}}$ over $\CC[u]$.

Each $\hc$ is a prime noetherian ring because its associated graded ring has these
properties \cite[Theorem 1.3]{EG}, \cite[Theorem 3.17]{Mont},
and \cite[Prop. 1.6.6, Theorem 1.6.9]{McCRob}.

\subsection{}
Let $\zc$ denote the centre of $\hc$, and set $\xc = \spec \zc$. It is known that $\zc$ is a finitely generated
$\CC$-algebra, and that $\hc$ is a finite $\zc$--module, \cite[Theorem 1.5 and Theorem 3.1]{EG}. By a lemma of
Dixmier, it follows that all simple $\hc$--modules are finite
dimensional. In fact $|\Gamma_n|$ is the strict upper bound for
the dimension of simple $\hc$--modules, \cite[Proposition 3.8]{EG}. We
therefore have a map
$$\chi : \textsf{Simp} (\hc) \longrightarrow \xc$$ which sends a
simple module $S$ to its central character $\chi(S)\in \xc$.

\subsection{}
The algebra $\zc$ has a Poisson bracket, making $\xc$ a Poisson
variety.
The following subsets of $X_{\mathbf{c}}$ are the same:
\begin{enumerate}
\item{}
the locus where the form is non-degenerate;
\item{}
its non-singular locus, $\textsf{Sm}(\xc)$ \cite[Theorem 7.8]{BG};
\item{}
the Azumaya locus of $\hc$, that is $\{\chi(S) : S$ is a simple
module of maximal dimension$\}$ \cite[Theorem 1.7]{EG};
\item{}
$\{x \in \xc : \chi^{-1}(x) \text{ is a singleton}\}$ \cite[Theorem 3.7] {EG}.
\end{enumerate}
Moreover, by \cite[Theorem 3.7]{EG}, if $x\in \textsf{Sm}(\xc)$,
then the unique simple $\hc$-module having central character $x$
is isomorphic to $\CC\Gamma_n$ as a $\CC\Gamma_n$-module.

\subsection{}
Let $e\in \CC\Gamma_n$ be the symmetrising
idempotent $|\Gamma_n|^{-1}\sum_{\gamma\in\Gamma_n} \gamma$.
The map $\zc \longrightarrow e\hc e$, $z\mapsto ze$, is an
isomorphism, \cite[Theorem 3.1]{EG}.
Thus, following \ref{ff}, there is a flat family of
commutative algebras $Z_{u\mathbf{c}}$ over $\CC[u]$. We set
$X_{u\mathbf{c}}= \spec Z_{u\mathbf{c}}$.

\subsection{}
\label{az}
For later use we need the following lemma.
\begin{lemma} For non--zero $\mathbf{c}$, the following sets are the same
\begin{enumerate}
\item the smooth locus, $\textsf{Sm}(X_{u\mathbf{c}})$;
\item the Azumaya locus of $H_{u\mathbf{c}}$;
\item the central characters of simple $\hc$-modules that are isomorphic to the
regular representation of $\Gamma_n$ as $\Gamma_n$-modules.
\end{enumerate}
\end{lemma}
\begin{proof}
As the generic simple module for $H_{u\mathbf{c}}$ has dimension
$|\Gamma_n|$ it follows that the Azumaya locus of
$H_{u\mathbf{c}}$ is the union of the Azumaya loci for
$H_{\lambda\mathbf{c}}$ for $\lambda \in \CC$. This proves the
equivalence of (1) and (3).

Since the associated graded ring of $H_{u\mathbf{c}}$ is
$\CC[V\oplus \CC]\ast \Gamma_n$, it follows from \cite[Proposition
1.6.6, Theorem 1.6.9 and Corollary 7.6.18]{McCRob} and
\cite[Theorem 3.8]{BrGoo} that the equivalence of (1) and (2)
follows if the non--Azumaya locus has codimension at least 2.
This, however, is clear as for each $\lambda\in \CC$ the Azumaya locus
of $H_{\lambda\mathbf{c}}$ has codimension at least 2 in
$X_{\lambda\mathbf{c}}$.
\end{proof}

\subsection{}
\label{norm} Observe as a consequence of the above proof that
$Z_{u\mathbf{c}}$ is normal. Indeed, $Z_{u\mathbf{c}}$ is Cohen--Macaulay (in
fact Gorenstein) since $eH_{u\mathbf{c}}e$ has an associated graded ring
isomorphic to $\CC [V \oplus \CC]^{\Gamma_n}$, \cite[Theorem 3.3]{EG}.
Combining this with the smoothness of $X_{u\mathbf{c}}$ in codimension one,
which is proved above, shows that $Z_{u\mathbf{c}}$ is normal, \cite[Theorem
2.2.11]{CG}.

\section{Non--commutative crepant resolutions}
\label{nc}
\subsection{}
Let $X$ be a Gorenstein $\CC$--variety.
A resolution of singularities $f: Y\longrightarrow
X$ is said to be \textit{crepant} if $f^*\omega_X = \omega_Y$.
Crepant resolutions are a generalisation of the notion of a
minimal resolution in two dimensions. However, crepant
resolutions need not exist, and need not be unique when they do
exist.

\subsection{}
\label{sympcrep} In the setup of Section \ref{sra} we could take
$X=\xc$. This is a Gorenstein variety, \cite[Theorem 1.5(i)]{EG}.
Moreover, since the Poisson form on $\xc$ is symplectic on
$\textsf{Sm}(\xc)$, the canonical bundle $\omega_{\xc}$ is
trivial. Thus any crepant resolution of $\xc$ has trivial
canonical class.
\begin{remark}
In fact, $\yc$ is a crepant resolution of $\xc$ if and only if
$\yc$ is a symplectic variety whose form agrees with that of $\xc$
on $\textsf{Sm}(\xc)$, \cite[Proposition 3.2]{kal0}
\end{remark}

\subsection{}
For the following definitions see \cite[Sections 3 and 4]{vdb2}. Throughout $R$ denotes a commutative noetherian domain over $\CC$. A module--finite $R$--algebra $A$ is
\textit{homologically homogeneous} if, for all $\mathfrak{p} \in \spec R$,
$\gldim A_{\mathfrak{p}} = \Kdim R_{\mathfrak{p}}$ and
$A_{\mathfrak{p}}$ is maximal Cohen--Macaulay.
A \textit{non--commutative crepant resolution} of $R$ is
a homologically homogeneous $R$--algebra of the form $A =
\edo_R(M)$, where $M$ is a reflexive $R$--module.

A justification for this definition of non--commutative crepant resolution is given in \cite[Section 4]{vdb2}.

\subsection{}
\label{stab}
We recall some material from \cite[Sections 3 and 4]{K} and
\cite[Section 6]{vdb2}. Set $X= \spec R$.
Let $A$
be an $R$--algebra that is finitely generated as an $R$--module, and let
$(e_i)_{i=1, \ldots ,p}$ be pairwise orthogonal idempotents in
$A$ such that $1= \sum_i e_i$.

For a map $R\longrightarrow K$ with $K$ a field and $V$ a finite
dimensional $A\otimes_R K$--module, we write $\underline{\dim} V =
(\dim_K e_iV)_i \in \mathbb{Z}^p$.

Pick $\lambda \in \homo_{\mathbb{Z}}(\mathbb{Z}^p,\mathbb{R})$ and let $\alpha = \underline{\dim}
V$.
We say that a finite dimensional $A \otimes_R K$-module
$V$ is \textit{stable} (respectively, \textit{semi-stable})
with respect to $\lambda$ if $\lambda (\alpha) = 0$ and
for every proper $A \otimes_RK$-submodule
$W$ of $V$ we have $\lambda (\underline{\dim} W) < 0$
(resp., $\lambda(\underline{\dim} W) \le 0$).

Our definition of stability is different from that in \cite{vdb2}:
where we have $\lambda(\underline{\dim} W)<0$ Van den Bergh has
$\lambda(\underline{\dim} W) >0$. Of course, the difference is
only cosmetic because we can pass back and forth between the two
notions by replacing $\lambda$ by $-\lambda$. The reason for this
difference is that later on in Section \ref{cluster} we want our
notion of stability to coincide with the notion that Haiman uses
in \cite{hai}.

We say that $\lambda$ is \textit{generic} (for $\alpha)$ if all semi-stable
representations of dimension vector $\alpha$ are stable.
There is a generic $\lambda$ if and only if $\alpha$ is indivisible, meaning that the
greatest common divisor of the $\alpha_i$s is 1.
The condition $\lambda  (\beta) \neq 0$ for all $0 < \beta
< \alpha$ ensures $\lambda$ is generic.

\subsection{}

Let $T$ be an $R$-scheme.
A \textit{family} of $A$-modules
of dimension $\alpha$
parametrised by $T$
is a locally free sheaf ${\mathcal F}$ of
${\mathcal O}_T$-modules
together with an $R$-algebra homomorphism $\phi:A \to \End_T
{\mathcal F}$
such that $e_i{\mathcal F}$ has constant rank $\alpha_i$ for all $i$.
We say that ${\mathcal F}$ is
{\it semi-stable} (resp., {\it stable}) if for every every field $K$
and every morphism $\xi:\Spec K \to T$,  $\xi^*{\mathcal F}$
is a semi-stable (resp., stable) $A \otimes_R K$-module.
Two families $({\mathcal F},\phi)$ and $({\mathcal F}',\phi')$
are \textit{equivalent} if there is an invertible
${\mathcal O}_T$-module ${\mathcal L}$ and an isomorphism $\psi:
{\mathcal F} \to {\mathcal F}' \otimes {\mathcal L}$ such that
the diagram
$$
\begin{CD}
A @>{\phi}>> \End {\mathcal F}
\\
@V{\phi'}VV @VV{\psiol}V
\\
\End {\mathcal F}' @>{\sim}>{\id \otimes {\mathcal L}}>
\End {\mathcal F}' \otimes {\mathcal L}
\end{CD}
$$
commutes; in this diagram $\psiol(\nu)=\psi\nu\psi^{-1}$.

A family $({\mathcal U},\rho)$ parametrised by $W$ is
{\it universal}
if for every $R$-scheme $T$ and every family $({\mathcal F},\phi)$
over $T$ there is a unique morphism $\xi:T \to W$ such that
$(\xi^*{\mathcal U},\xi^* \rho)$ is equivalent to
$({\mathcal F},\phi)$; here $\xi^*  \rho$ denotes
the composition $A \to \End {\mathcal U} \to \End \xi^*
{\mathcal U}$.

Proposition \ref{prop.vdb} below says that in suitable situations
there is a universal family. We also express this by saying that $W$ is
a fine moduli space for families of $A$-modules of dimension $\alpha$
and call ${\mathcal U}$ the {\it universal family}.

\subsection{}
\label{repres}

%Fix $\alpha\in \mathbb{Z}^p$ and suppose that
%$\lambda$ is generic for $\alpha$.

Define a functor
$\mathbf{R}^s: \hbox{$R$-schemes} \longrightarrow \Sets$ by
$$
T \mapsto \{\hbox{equiv. classes of families of
$\lambda$--stable $A$--modules over $T$ with dimension $\alpha$}\}.
$$

\begin{prop}
\label{prop.vdb}
\cite[Proposition 6.2.1]{vdb2}
Suppose that $\lambda$ is generic for $\alpha$. Then {\bf R}${}^s$
is represented by a closed subscheme $W^s \subset \PP_X^N$.
\end{prop}

We write $f:W^s \to X$ for the structure morphism, and
${\mathcal B}$ for the universal family
of $\lambda$-stable $A$-modules of dimension $\alpha$.

If $U$ is either an open or closed subscheme of $X$ then the
representing scheme for $U$ is $f^{-1}(U)$ and its universal family is
${\mathcal B} \big\vert{}_U$ (cf. the sentence after \cite[Lemma 6.2.2]{vdb2}).

\subsection{} \label{trivspace}
It is shown in \cite[Lemma 6.2.3]{vdb2} that in the
case $A=\mat_{\sum\alpha_i}(R)$ the map $W^s\longrightarrow X$ is
an isomorphism.

\subsection{}
\label{der}

Assume $A=\End_R M$ is non--commutative crepant resolution of $R$.
Suppose we have an $R$-module decomposition $M = \oplus_{i=1}^p M_i$ and
let $e_1,\ldots,e_p$ be the projections onto the $M_i$s viewed as
idempotents in $A$.
Define
$$
\alpha_i:=\rank M_i=\rank e_iM.
$$
Suppose that $\lambda$ is generic for $\alpha$. Let
$f: W^s\longrightarrow X$ and $\mathcal{B}$ be as in \ref{repres}.

Let $U \subset X$ be the locus where $M$ is locally free. It
follows from \ref{trivspace} that $f^{-1}(U)\longrightarrow U$
is an isomorphism. Let $Y$ be the closure of $f^{-1}(U)$; this is
the unique irreducible component of $W^s$ mapping birationally onto $X$.
We continue to denote the restriction of $f$ to
$Y$ by $f$. Let $\mathcal{G}$ be the restriction of $\mathcal{B}$
to $Y$.

There is a pair of adjoint functors between $D^b(\coh Y)$ and
$D^b(\md A)$:
\begin{eqnarray*}
\Phi: D^b(\coh Y)\longrightarrow D^b(\md A) : C \mapsto
\mathbf{R}\Gamma(C\otimes^{\mathbf{L}}_{\mathcal{O}_Y} \mathcal{G}) \\
\Psi : D^b(\md A) \longrightarrow D^b(\coh Y) : D \mapsto
D\otimes^{\mathbf{L}}_A\mathcal{G}^*.
\end{eqnarray*}
The following theorem follows \cite{BKR} closely.
\begin{theorem}\cite[6.3.1]{vdb2}
Assume that for every point $x\in X$,
$$
\dim (Y\times_X Y)\times_X \spec \mathcal{O}_{X,x}  \le \codim x +1.
$$
Then $f: Y \longrightarrow X$ is
a crepant resolution of $X$ and $\Phi$ and $\Psi$ are inverse
equivalences.
\end{theorem}

\subsection{}
\label{fibreres}
Let $x\in X$.
We write $D^b_x(\coh Y)$ and $D^b_x(\md A)$ for the full
subcategories of $D^b(\coh Y)$ and $D^b(\md A)$ consisting of
complexes supported on $f^{-1}(x)$ and on $x$ respectively.
Since $\Phi$ and $\Psi$ are functors over $X$ they restrict to
equivalences between $D^b_x(\coh Y)$ and $D^b_x(\md A)$,
\cite[9.1]{BKR} and \cite[6.6]{vdb2}.

\subsection{}

Fix a non--zero $\mathbf{c}$.
The following lemma allows us to apply the machinery in this section.

\begin{lemma}
\label{lem.nc.crep}
The algebra $\hc$ is a non--commutative crepant resolution of $\xc$.
\end{lemma}
\begin{proof}
By \cite[Theorem 15.(ii), (iii)]{EG} $\hc e$ is a finitely generated, reflexive,
Cohen--Macaulay $\hc$--module. Furthermore, by \cite[Theorem 1.5(iv
)]{EG} we have $\hc \cong \edo_{e\hc e}(\hc e)$.
By \cite[Corollary 6.18]{McCRob}, $\gldim \hc < \infty$.
Thus $\hc$ is a non--commutative crepant resolution of
$e\hc e \cong \zc$ by \cite[Lemma 4.2]{vdb2}.
\end{proof}

Let $H_{u\mathbf{c}}$ be the flat family of symplectic reflection algebras
defined in \ref{ff}. The arguments used to prove Lemma \ref{lem.nc.crep}
all extend to $H_{u\mathbf{c}}$, showing that $H_{u\mathbf{c}}$
is a non-commutative crepant resolution of $X_{u\mathbf{c}}$.

\subsection{}
\label{sect.apply-to-Hc}

If we set $R=eH_{u\mathbf{c}}e$, $M=H_{u\mathbf{c}}e$, and
$A=H_{u\mathbf{c}}$, we are in the situation of subsection
\ref{der}. Let $(e_i)_{i=1, \ldots ,p}\in \CC\Gamma_n$ be the
central orthogonal idempotents corresponding to the irreducible
representations, labelled so that $e_1$ corresponds to the trivial
representation. Thus $e_1=e$. By the PBW theorem for $\hc$ and
$H_{u\mathbf{c}}$ we can consider the $e_i$s as elements of $\hc$
and $H_{u\mathbf{c}}$. If we set $M_i=e_iM$ there is a
decomposition $H_{u\mathbf{c}}e=M=\oplus_{i=1}^p M_i$ as in
Section 3.8. Let $\alpha_i=\rank e_i H_{uc}e$ and write
$\alpha=(\alpha_i)_{1 \le i \le p} \in \ZZ^p$.

Let $E_i$ be the simple $\CC\Gamma_n$-module corresponding to $e_i$.
As remarked in the proof of \cite[Lemma 2.24]{EG},
$e_iH_{u\mathbf{c}}e \cong
\Hom_{\Gamma_n}(E_i,H_{u\mathbf{c}}e)$
so
$$
\alpha_i=\rank_{Z_{u\mathbf{c}}} e_iH_{u\mathbf{c}}e = \dim E_i.
$$
In particular, $\alpha_1=1$ so $\alpha$ is indivisible, and there are many
maps $\lambda:\ZZ^p \to \RR$ such that $\lambda(\alpha)=0$ and
$\lambda(\beta)>0$ for all $ 0<\beta<\alpha$.
For each such $\lambda$ there is a moduli space for $\lambda$--stable
$H_{u\mathbf{c}}$--modules having dimension $\alpha$ (equivalently, that are
isomorphic to $\CC\Gamma_n$ as $\CC\Gamma_n$-modules), and
each such moduli space has a unique irreducible component that maps
birationally to $X_{u\mathbf{c}}$.

\section{Semi--small maps}
\label{ss}
\subsection{}
\label{defss}
A proper birational map $f:Y \longrightarrow X$ between irreducible
varieties is \textit{semi-small} if
$$
2 \codim_Y Z \geq \codim_X f(Z)
$$
for all irreducible subvarieties $Z \subset Y$.
Note that if $f$ is semi--small the above inequality holds for all (not necessarily irreducible) subvarieties of $Y$. It is a theorem of Verbitsky, \cite[Theorem 2.8]{ver}, and Kaledin, \cite[Proposition 4.4]{kal0}, that any crepant resolution of $V/\Gamma_n$ is semi--small.
\subsection{}
\label{sathyp}

The following result relates semi--small maps to the hypothesis
of Theorem \ref{der}.
\begin{lemma}
Suppose that $f: Y\longrightarrow X$ is a semi--small morphism
between irreducible varieties of finite type over a field $k$. Then
$$
\dim (Y\times_X Y)\times_X \spec \mathcal{O}_{X,x} \leq \codim x +1
$$
for every point $x\in X$.
\end{lemma}
\begin{proof}
It is well-known (see, e.g., \cite[Proposition 2.1.1 and Remark 2.1.2]{dCM})
that semi-smallness of $f$ is equivalent to the condition that every
irreducible component of $Y \times_X Y$ has dimension at most $\dim
X$, so it suffices to prove that if $Z$ is an irreducible variety
of dimension $\le \dim X$ and $g:Z \to X$ a morphism,
then $\dim Z \times_X \Spec {\mathcal O}_{X,x} \le \codim x$ for every
point $x \in X$.

This reduces to the affine case.
We need to prove the following: if $R \to S$ is a homomorphism
between two domains that are finitely generated $k$-algebras such that
$\Kdim S \le \Kdim R$, then
$\Kdim S \otimes_R R_\fp \le \Kdim R_\fp$ for all $\fp \in \Spec R$.

Let $\Rol$ denote the image of $R$ in $S$, and let $\fq \in \Spec R$
denote the kernel of $R \to S$.
The hypotheses on $R$ and $S$ are such that they, their
localizations, and the homomorphic images of these are catenary.

Write $d$ for $\Kdim S-\Kdim \Rol$.
Since the Krull dimensions of $\Rol$ and $S$ are equal to the
transcendence degrees over $k$ of their fraction fields,
we can write  $S=\Rol[x_1,\ldots,x_d,\ldots,x_n]$ where
$\{x_1,\ldots,x_d\}$ is algebraically independent over $\Rol$,
and $x_{d+1},\ldots,x_n$ are algebraic over $\Rol[x_1,\ldots,x_d]$.

Write $\Rol_{\fp}$ for $\Rol \otimes_R R_{\fp}$. Then
$S \otimes_R R_\fp = \Rol_{\fp}[x_1,\ldots,x_n]$. The
Krull dimension of an extension $C[x]$ is equal to
either $\Kdim C$ or $\Kdim C +1$
depending on whether $x$ is algebraic over $C$ or not,
so an induction argument
shows that $\Kdim S \otimes_R R_\fp \le \Kdim \Rol_{\fp} +d$.
Thus
\begin{align*}
\Kdim S \otimes_R R_\fp
& \le \Kdim \Rol_{\fp} + \Kdim S - \Kdim \Rol
\\
& \le \Kdim \Rol_{\fp} + \Kdim S- \Kdim R + \hgt \fq,
\end{align*}
where $\hgt \fq$ denotes the height of $\fq$.
But $\hgt \fq \le \hgt \fq R_{\fp}$ so
$$
\Kdim S \otimes_R R_\fp
%& \le \Kdim \Rol_{\fp} + \Kdim S- \Kdim R + \hgt \fq
%\\
 \le \Kdim \Rol_{\fp} + \Kdim S- \Kdim R + \hgt \fq R_{\fp}
= \Kdim R_{\fp} + \Kdim S- \Kdim R,
$$
and the result follows because $\Kdim S \le \Kdim R$.
\end{proof}

\subsection{}
\label{defdim}
For the rest of the section we will be interested in schemes of finite type over $\CC$ with a $\CC^*$-action. All morphisms will be
$\CC^*$-equivariant. We always consider $\CC$ as to have the multiplicative action of $\CC^*$.
In case $f: X\longrightarrow \CC$ is a $\CC^*$--equivariant morphism, we will denote the fibres $f^{-1}(s)$ by $X_s$. Note that $X_s \cong X_1$ for all non--zero $s\in \CC$.

We say that an affine variety $X$ has an \textit{expanding} $\CC^*$--action if the corresponding $\mathbb{Z}$--grading $\CC[X] = \oplus_{i\in\mathbb{Z}} \CC[X]_i$ is concentrated in non--negative degrees.

\subsection{} The following lemma can be compared with \cite[II.4.2, Satz 2]{kraft}.
\begin{lemma}
\label{lem.kraft}
Let $X$ be an irreducible variety and  $f:X \longrightarrow \CC$ be a $\CC^*$-equivariant morphism.
Let $Z \subset X_s$ be an irreducible subvariety. Then either $\overline{\CC^* Z}\cap X_0$ is empty or
$$
\dim
(\overline{\CC^* Z} \cap X_0)
=\dim Z.
$$

Furthermore, if $X$ is an affine variety with expanding $\CC^*$--action, then $\overline{\CC^* Z}\cap X_0$ is non--empty.

\end{lemma}
\begin{proof}
Suppose $\overline{\CC^* Z}\cap X_0$ is non--empty. Since $\overline{\CC^*Z}$ is irreducible and  $\dim (\overline{\CC^*Z}) = \dim Z +1$, we have $\dim(\overline{\CC^*Z}\cap X_0)\leq \dim Z$. On the other hand,
the dimension of the fibers of the restriction $f|_{\overline{\CC^*Z}}:
\overline{\CC^* Z} \to \CC$ is minimal on a dense open
set of $\mathbb{C}$. Since the $\CC^*$--action identifies the fibres of this map over non--zero elements of $\CC$, we see that the minimal fibre dimension is bounded below by $\dim Z$, as required.

Now assume that $X$ is an affine variety with expanding $\CC^*$--action. Let $I$ be the ideal of $\mathbb{C}[X]$ annihilating $Z$. Then the ideal corresponding to $\overline{\CC^*Z}\cap X_0$ is $\gr I$, the ideal consisting of leading terms of elements of $I$, \cite[II.4.2, Satz 3]{kraft}. In particular, as $I$ is proper, so too is $\gr I$. Thus $\overline{\CC^*Z}\cap X_0$ is non--empty.
\end{proof}

\subsection{}
We need a simple lemma.
\label{push}
\begin{lemma}
\label{lem.good.pi}
Suppose there is a commutative diagram of $\CC^*$-equivariant morphisms
$$
\setlength{\unitlength}{1mm}
\begin{picture}(25,15)
\put(0,13){$Y$}
\put(4,15){\vector(1,0){20}}
\put(25,13){$X$}
\put(14,16){$\scriptstyle\pi$}
\put(3,11){\vector(1,-1){7}}
\put(24,11){\vector(-1,-1){7}}
\put(12,0){$\CC$}
\end{picture}
$$
where $\pi$ is proper.
If $Z \subset Y_s$ is an irreducible subvariety, then
\begin{equation}
\label{eqn.good.pi}
\pi(\overline{\CC^* Z} \cap Y_0) =
\overline{\CC^* \pi(Z)} \cap X_0.
\end{equation}
In particular, if $\overline{\CC^* Z} \cap Y_0$ is non--empty, then $\dim (\overline{\CC^* Z} \cap Y_0) = \dim Z$ and $\dim \pi(\overline{\CC^* Z} \cap Y_0) = \dim \pi(Z)$.
\end{lemma}
\begin{proof}
First we show that the right-hand side of (\ref{eqn.good.pi}) is contained
in the left-hand side.
Let $x \in \overline{\CC^* \pi(Z)} \cap X_0$. Certainly, $\CC^*\pi(Z)=\pi (\CC^* Z)
\subset \pi(\overline{\CC^* Z})$; but the last term is closed because $\pi$
is proper, so $\overline{\pi(\CC^*Z)} \subset \pi(\overline{\CC^* Z})$.
Hence $x =\pi(y)$ for some $y \in \overline{\CC^* Z}$; but $x \in X_0$, so $y \in
Y_0$, whence $x \in \pi(\overline{\CC^* Z} \cap Y_0)$.

The proof that the left-hand side of (\ref{eqn.good.pi}) is contained
in the right-hand side does not depend on the hypothesis that $\pi$
is proper:
For any subsets $W$ and $W'$ of $Y$, $\pi(W \cap W') \subset \pi(W) \cap \pi(W')$,
and $\pi(\overline{W}) \subset \overline{\pi(W)}$. Thus $\pi(\overline{W} \cap W') \subset \pi(\overline{W})
\cap \pi(W') \subset\overline{\pi(W)} \cap \pi (W')$. Now apply this with
$W=\CC^* Z$ and $W'=Y_0$.

Under the non--emptiness hypothesis, the equality of dimensions follows from Lemma \ref{defdim}.
\end{proof}

\subsection{}
The following result allows us to deform semi--small morphisms.
\begin{lemma}
\label{lem.good.pi.2}
Let $X$ be an affine variety with expanding $\CC^*$--action and suppose that we have a commutative diagram of $\CC^*$-equivariant morphisms
$$
\setlength{\unitlength}{1mm}
\begin{picture}(25,15)
\put(0,13){$Y$}
\put(4,15){\vector(1,0){20}}
\put(25,13){$X$}
\put(14,16){$\scriptstyle\pi$}
\put(3,11){\vector(1,-1){7}}
\put(24,11){\vector(-1,-1){7}}
\put(12,0){$\CC$}
\end{picture}
$$
where $\pi$ is proper. Assume that $\dim Y_0 = \dim Y_s$ and $\dim X_0 = \dim X_s$ for all $s\in\CC$. If $\pi_0: Y_0 \longrightarrow X_0$ is semismall, so too is $\pi_s:Y_s \longrightarrow X_s$ for
all $s$.
\end{lemma}
\begin{proof}
Let $Z \subset Y_s$ be an irreducible subvariety. Set $Z_0:=\overline{\CC^* Z} \cap Y_0$. As $X$ has an expanding action, $\pi(Z_0) = \overline{\CC^*\pi(Z)}\cap X_0$ is non--empty by Lemmas \ref{defdim} and \ref{push}. Thus the semi--small hypothesis shows
\begin{equation}
\label{eqn.1/2-small}
2 \dim Y_0 - 2 \dim Z_0 \ge \dim X_0 - \dim \pi(Z_0).
\end{equation}
By Lemma \ref{push}, we may replace $\dim Z_0$ and $\dim \pi(Z_0)$ in this
inequality by $\dim Z$ and $\dim \pi(Z)$. Now the lemma follows by replacing $\dim Y_0$ and $\dim X_0$ in this inequality by $\dim Y_s$ and $\dim X_s$.
\end{proof}

\section{Application}
\label{main}
\subsection{}
\label{cluster}
A representation of $H_{\mathbf{0}} = \CC[V]\ast \Gamma_n$ is
called a $\Gamma_n$--\textit{constellation} if its restriction to
$\CC \Gamma_n$ is isomorphic to the regular representation.
A constellation $M$ is a \textit{cluster} if it is generated as a
$\CC[V]$ module by $M^{\Gamma_n}$, the copy of the
trivial representation it contains.

We write $K(\Gamma_n)$ for the Grothendieck group of $\CC\Gamma_n$ and
$\alpha$ for the class of regular representation.
Let $\lambda : K(\Gamma_n) \longrightarrow \mathbb{R}$ be a linear
function such that $\lambda(\alpha)=0$.
Following \ref{stab}, a constellation $M$ is $\lambda$
(semi--)stable if for every proper $H_{\mathbf{0}}$--submodule
$N\subset M$ we have $\lambda (N) (\leq ) <0$.
For generic $\lambda$ there is a moduli space of
$\lambda$--stable $\Gamma_n$--constellations, a projective scheme
over $\spec Z_{\mathbf{0}} = V/\Gamma_n$.

If $\lambda$ is chosen so that $\lambda(\alpha)=0$ and, for each simple
$\Gamma_n$-module $S$,
$$
\lambda(S)=\begin{cases}
    1 & \text{if $S$ is trivial}
    \\
    < 0 & \text{if $S$ is not trivial}
\end{cases}
$$
then a constellation is $\lambda$-stable if and only if it is
a cluster.

\subsection{}
\label{haith} The following construction was given in
\cite[Corollaries 3 and 4]{wang2}. Let $X_{\Gamma}$ be the minimal
resolution of the Kleinian singularity $\CC^2/\Gamma$. We have
maps
\begin{eqnarray*} \hil^n(X_{\Gamma}) \longrightarrow \sym^n(X_{\Gamma}) \longrightarrow \sym^n(\CC^2/\Gamma) \cong V/\Gamma_n, \end{eqnarray*}
where the first map is the Hilbert--Chow map, \cite[Chapter
1]{nakbook}, and the second arises from functoriality. Since
$X_{\Gamma}$ is symplectic, so too is $\hil^n(X_{\Gamma})$ and
thus the composition is a crepant resolution of $V/\Gamma_n$, see
\ref{sympcrep}.

\subsection{}
Let $\hil^{S_n}(X_{\Gamma}^n)$ denote the $S_n$--Hilbert scheme of
Ito and Nakamura, \cite[Introduction and Sect. 8.2]{in}. The
following result is due to Haiman: we include our own outline of
the proof for the reader's benefit.

\begin{theorem}\cite[Section 7.2.3]{hai2}, \cite{hai3} There is an isomorphism between
$\hil^n(X_{\Gamma})$ and $\hil^{S_n}(X_{\Gamma}^n)$. In particular,
there exists $\lambda$ such that $\hil^n(X_{\Gamma})$ is a moduli space of $\lambda$--stable $\Gamma_n$ constellations.
\end{theorem}
\begin{proof} The isomorphism follows from the $n!$--conjecture applied to
the smooth surface $X_{\Gamma}$, \cite[Sect. 5.2]{hai} and \cite{hai2}.
There is a  commutative diagram of $S_n$-equivariant morphisms
\[
\begin{CD} \Sigma @> p >> X_{\Gamma}^n \\
@V q VV @VVV \\
\hil^{S_n}(X^n_{\Gamma}) @>>> V/\Gamma_n
\end{CD}
\]
in which $\Sigma$ is the reduced fibre product and
$q$ is finite and flat.

It is well-known that $X_{\Gamma}$ is a fine moduli space for
$\Gamma$-clusters on $\CC^2$. We write ${\mathcal B}$ for the
locally free sheaf on $X_\Gamma$ that is the
universal family of $\Gamma$-clusters.
The obvious permutation action makes ${\mathcal B}^{\boxtimes n}$
an $S_n$-equivariant sheaf on $X_{\Gamma}^n$.
Let $\mathcal{P}$ denote
%the bundle on $\hil^n (X_{\Gamma})$ corresponding to
$q_*p^*\mathcal{B}^{\boxtimes n}$.
Because $q$ and $p$ are $S_n$-equivariant ${\mathcal P}$ is an
$S_n$-equivariant sheaf on $\hil^{S_n}(X^n_{\Gamma})$. Since
the $S_n$-action on  $\hil^{S_n}(X^n_{\Gamma})$ is trivial this means
that $S_n$ acts as automorphisms of ${\mathcal P}$.

The ring homomorphism $\CC[x,y]*\Gamma \to \End {\mathcal B}$
induces homomorphisms $\CC[x,y]^{\otimes n} * \Gamma^n \to \End
p^* {\mathcal B}^{\boxtimes n}$ and $\CC[V] * \Gamma^n \to \End
q_*p^* {\mathcal B}^{\boxtimes n} = \End {\mathcal P}$. Combining
the last of these with the $S_n$-action produces a ring
homomorphism $\CC[V]*\Gamma_n \to \End{\mathcal P}$.

Consider $\lambda:K(\Gamma_n) \to \RR$ of the form
$$\lambda (M) = C \rho (M|_{\Gamma^n}) + \sigma (M|_{S_n})$$
where $\rho:K(\Gamma^n) \to \RR$ and $\sigma:K(S_n) \to \RR$ are
such that stable constellations are clusters. It can be shown that
for a suitable choice of $C \gg 1$, the geometric fibers
$\mathcal{P}(x):=\mathcal{P}/\fm_x\mathcal{P}$ of ${\mathcal P}$
are $\lambda$-stable $\Gamma_n$--constellations, and hence  that
$\mathcal{P}$ is a family of $\lambda$--stable
$\Gamma_n$--constellations.

The fixed points subsheaf $\mathcal{P}^{\Gamma_{n-1}}$ is the
universal family of $\CC[x,y]*\Gamma$-modules whose fibres have
$n$ copies of the regular representation of $\Gamma$ \cite[Prop.
7.2.12]{hai2}.

Let $M_{\lambda}$ be the moduli space of $\lambda$--stable
$\Gamma_n$--constellations and ${\mathcal S}$ the
universal family on it of $\lambda$-stable
$\CC[V]*\Gamma_n$-constellations.

The homomorphisms $\phi: \CC[V]*\Gamma_n \to \End {\mathcal P}$
and $\psi:\CC[V]*\Gamma_n \to \End {\mathcal S}$ restrict to
homomorphisms $\phi':\CC[x,y]*\Gamma \to \End {\mathcal
P}^{\Gamma_{n-1}}$ and $\psi':\CC[x,y]*\Gamma \to \End {\mathcal
S}^{\Gamma_{n-1}}$.

Since $\mathcal{P}$ is a family of $\lambda$--stable
$\Gamma_n$--constellations, there is a morphism
$f:  \hil^n(X_{\Gamma}) \longrightarrow M_{\lambda}$
such that $f^*\mathcal{S} \cong \mathcal{P}$ and $\phi=f^* \psi$.
Thus $\phi'=f^*  \psi'$.
Similarly, by the universal property of
${\mathcal P}^{\Gamma_{n-1}}$, there exists
$g: M_{\lambda} \longrightarrow \hil^n(X_{\Gamma})$
such that $g^*\mathcal{P}^{\Gamma_{n-1}} \cong \mathcal{S}^{\Gamma_{n-1}}$
and $\phi'=g^* \psi'$.

Both $M_{\lambda}$ and $\hil^n(X_{\Gamma})$ have trivial $\Gamma_n$-action so
$f$ and $g$ are automatically $\Gamma_n$--equivariant. Therefore
$$
(gf)^*\mathcal{P}^{\Gamma_{n-1}} \cong f^*(\mathcal{S}^{\Gamma_{n-1}})
\cong \mathcal{P}^{\Gamma_{n-1}}.
$$
Notice too that $f^*g^*  \psi'=\psi'$.
Since $\hil^n(X_{\Gamma})$ is a fine moduli space with universal
family $\mathcal{P}^{\Gamma_{n-1}}$, it follows that $gf=\id$.

There is a non-empty open subset $U$ of $V/\Gamma_n$ such that the natural maps
$\alpha:\hil^n X_\Gamma \to V/\Gamma_n$ and
$\beta:M_\lambda \to V/\Gamma_n$ restrict to
isomorphisms $\alpha^{-1}(U) \to U$ and $\beta^{-1}(U) \to U$.
The closure $Y_\lambda$ of $\beta^{-1}(U)$ is the unique irreducible
component of $M_\lambda$ that maps birationally to $V/\Gamma_n$.
Since $f$ and $g$ are morphisms of $V/\Gamma_n$-schemes and $fg=\id$ they
restrict to mutually inverse isomorphisms between $\alpha^{-1}(U)$
and $\beta^{-1}(U)$
and hence between their closures. But $\hil^n X_\Gamma$ is irreducible, so
$f$ and $g$ yield an isomorphism $\hil^n X_\Gamma \cong Y_\lambda$.
\end{proof}

\subsection{}
\label{haith2}

As noted in \cite[Section 4.4]{wang1} and \cite{hai2} the previous
theorem, together with the main result in \cite{BKR}, has the
following important consequence.

\begin{corollary}
\label{cor.haiman}
The derived categories $D^b(\coh \hil^n(X_{\Gamma}))$ and
$D^b(\md H_{\mathbf{0}})$ are equivalent.
\end{corollary}
\begin{proof}
This follows from Theorem \ref{haith} since the resolution
$\hil^n(X_{\Gamma})\longrightarrow V/\Gamma_n$ is crepant,
hence semismall by \ref{defss}, and so, using Lemma \ref{sathyp},
satisfies the hypothesis of Theorem \ref{der}.
\end{proof}

\subsection{}
\label{sect.perturb}
Set $\Theta = \{ \mu :K(\Gamma_n) \longrightarrow \mathbb{R} :
\mu (\alpha) = 0 \}$. Let $\lambda$ be the element in $\Theta$ given
by Theorem \ref{haith}. Define $\Theta_{\lambda}^+ = \{ 0\neq
M\subset \CC\Gamma_n : \lambda(M) >0 \}$ and
$\Theta_{\lambda}^- = \{ 0 \neq
M \subset \CC\Gamma_n : \lambda(M)<0\}$. If $M$ is a proper
$\Gamma_n$--submodule of the regular representation such that
$\lambda (M) = 0$ we can perturb $\lambda$ to $\mu$ so that
$\Theta_{\lambda}^+\cup \{ M\} \subseteq \Theta_{\mu}^+ $ and
$\Theta_{\lambda}^- \subseteq \Theta_{\mu}^-$.
The $\lambda$--stable constellations are the same as
the $\mu$--stable constellations
since $\lambda$ is generic.
Notice that every $\mu$-semistable constellation is stable.
Thus, without loss of generality, we may replace $\lambda$ by $\mu$
and assume that $\lambda(M) \neq 0$ for all proper
subrepresentations of a $\lambda$-stable constellation $M$.

\subsection{}

Let $H_{u\mathbf{c}}$ be the flat family of symplectic reflection
algebras defined in \ref{ff}.

We will now use Van den Bergh's result in Theorem \ref{der} to
extend Corollary \ref{cor.haiman} to the deformations $Y_{\mathbf{c}}
\to X_{\mathbf{c}}$ where $Y_{\mathbf{c}}$ is a suitable moduli space of
$H_{\mathbf{c}}$-modules.

In Section \ref{sect.apply-to-Hc}
 $\alpha$ denoted the element of $\ZZ^p$ defined by
$$
\alpha_i=\rank_{Z_{u\mathbf{c}}} e_iH_{u\mathbf{c}}e = \dim E_i
$$
where $E_i$ is the irreducible representation of $\Gamma_n$
corresponding to the central idempotent $e_i \in \CC\Gamma_n
\subset H_{\mathbf{c}}$. Therefore under the isomorphism
$K(\Gamma_n) \to \ZZ^p$, $[E_i] \mapsto
(\delta_{1i},\ldots,\delta_{ni})$, we have $[\CC\Gamma_n] \mapsto
\alpha$. In particular, the use of $\alpha$ in the last few
subsections is compatible with the use of $\alpha$ in Section
\ref{sect.apply-to-Hc}.

Let $\lambda:K(\Gamma_n) \to \RR$ be generic for $\alpha$. Let $W$
be the moduli space, as constructed in Section \ref{nc}, of
$\lambda$--stable $H_{u\mathbf{c}}$--modules isomorphic to
$\CC\Gamma_n$ (equivalently, of dimension $\alpha$), and let $Y$
be the irreducible component of $W$ that maps birationally to
$\textsf{Sm}(X_{u\mathbf{c}})$.

There is a natural $\CC^*$-equivariant map $W \to \CC$ and its restriction to
$Y$ fits into the following commutative diagram in the $\mathbb{C}^*$--equivariant
category
$$
\setlength{\unitlength}{1mm}
\begin{picture}(25,15)
\put(0,13){$Y$}
\put(4,15){\vector(1,0){20}}
\put(25,13){$X$}
\put(14,16){$\scriptstyle\pi$}
\put(3,11){\vector(1,-1){7}}
\put(2,6){$f$}
\put(24,11){\vector(-1,-1){7}}
\put(21,6){$g$}
\put(12,0){$\CC$}
\end{picture}
$$
The horizontal arrow is obtained by taking the central character.

\begin{theorem}
Keep the above notation.
\begin{enumerate}
\item The fibre $Y_{\mathbf{c}} := f^{-1}(1)$ is a crepant resolution of $\xc = \spec \zc$.
\item There is an equivalence of categories between $D^b (\coh Y_{\mathbf{c}})$ and $D^b(\md \hc)$.
\end{enumerate}
\end{theorem}
\begin{proof}
For $\tau\in \CC$, we write $W_{\tau}$ for the fiber of $W \to \CC$ over $\tau$;
we also define $Y_{\tau} = f^{-1}(\tau)$ and $X_{\tau} = g^{-1}(\tau)$.
By Theorem \ref{der} and Lemma \ref{sathyp}, it is enough to show that $Y_1$
is the irreducible component of $W_1$ that is birational to $\textsf{Sm}(X_1)$,
and that the restriction of $\pi$ to a morphism $Y_1\longrightarrow X_1$ is
semi--small.

The variety $Y_{\tau}$ is a moduli space of $\lambda$--stable
$H_{\tau\mathbf{c}}$--modules of dimension
$\alpha=\underline{\dim} (\CC\Gamma_n)$.
By Lemma \ref{az}, $Y_{\tau}$ contains the irreducible component of
$W_{\tau}$ that maps birationally to $\textsf{Sm}{X_{\tau}}$. In particular,
$\dim Y_{\tau} \geq  \dim X_{\tau}$; in fact, these dimensions
are equal because $Y$ is irreducible of dimension $\dim X = \dim X_{\tau}+1$.
If $\tau$ is non--zero, then $Y_{\tau}$ is irreducible since
$\overline{\CC^*Y_{\tau}}$ is a closed subset of Y of the same dimension as
$Y$. On the other hand, if $\tau = 0$ then \cite[Section 8]{BKR} combined with Corollary \ref{haith2} shows that this particular irreducible component is a connected component of the moduli space $W_0$.
By \ref{norm}, $X$ is normal, so Zariski's main theorem implies that $Y_{0}$
is connected \cite[Corollary III.11.4]{Ha}, and we see that $Y_{0}$ is also
irreducible. We deduce that for any $\tau$, $Y_\tau$ is the irreducible
component of $W_\tau$ that is birational to $\textsf{Sm}(X_{\tau})$.

The semi--smallness follows from Lemma \ref{lem.good.pi.2} since the restriction of $\pi$ to $Y_{0}\longrightarrow X_{0}$ is semi--small, being the crepant resolution of Theorem \ref{haith}.
\end{proof}

\subsection{}
\label{nosimples}
Let $\mathfrak{m}_x$ be the maximal ideal of $\zc$ corresponding to $x\in \xc$.
The simple $\hc$-modules with central character $x$ are precisely the
simple modules of the finite dimensional algebra $\hc/\mathfrak{m}_x\hc$.
\begin{corollary}
Let $\pi_{\mathbf{c}} : Y_{\mathbf{c}} \longrightarrow \xc$ be the crepant resolution above. There is an equivalence of
triangulated categories between $D^b_x(\coh Y_{\mathbf{c}})$ and $D^b_x(\md \hc)$. In particular there is an isomorphism between the Grothendieck groups $K(\pi_{\mathbf{c}}^{-1}(x))$ and $K(\hc/\mathfrak{m}_x\hc)$.
\end{corollary}
\begin{proof}
The first sentence has already been noted in \ref{fibreres}. By
devissage, the Grothendieck groups of $D^b_x(\coh Y_{\mathbf{c}})$
and $D^b_x(\md \hc)$ are isomorphic to
$K(\pi_{\mathbf{c}}^{-1}(x))$ and $K(\hc/\mathfrak{m}_x\hc)$
respectively, thus confirming the second sentence.
\end{proof}

\subsection{}
When $n=1$ the varieties $\xc$ are deformations of Kleinian singularities
$\CC^2/\Gamma$ and any crepant resolution coincides with the minimal resolution.
Hence, by the McKay correspondence, the $K$--theory of the fibre $f^{-1}(x)$ is completely determined by the type of orbifold singularity at $x\in \xc$. Indeed, suppose the singularity at $x\in \xc$ is locally of the form $\CC^2/G$ for some finite subgroup $G$ of $SL_2(\CC)$. Then the rank of $K(f^{-1}(x))$ equals the number of irreducible representations of $G$. On the other hand, the algebras $\hc$ are deformed preprojective algebras. These algebras also depend only on the type of orbifold singularity at $x\in \xc$, since there is a ''slice" theorem which reduces the representation theory to the case of the point $0\in \CC^2/G$, \cite[Corollary 4.10]{CB}. Thus the rank of $K(\hc /\mathfrak{m}_x\hc)$ also equals the number of irreducible representations of $G$, as expected.

Of course, in the case of arbitrary $n$ but $\mathbf{c}=0$, the corollary (which is an immediate consequence of Haiman's work) recovers the generalised McKay correspondence proved by Kaledin, \cite{kal}, for the orbifold singularities appearing locally in $V/\Gamma_n$.

\subsection{}
Corollary \ref{nosimples} gives us a geometric description for the number of simple modules in $\hc/\mathfrak{m}_x\hc$. In practice this is not immediately applicable as we have no geometric understanding of $Y_{\mathbf{c}}$. There is, however, evidence to suggest that the following question has a positive answer:

\begin{question} Let $V$ have Gorenstein singularities and let $v\in V$. Is $\rank K(f^{-1}(v))$ independent of the choice of crepant resolution $f: \tilde{V}\longrightarrow V$?
\end{question}

Indeed, \cite[Proposition 6.3.2]{denloe} shows that the mixed Hodge polynomial (of Borel--Moore homology) of $f^{-1}(x)$ is independent of the choice of resolution. Thus if the homology groups of the fibres are spanned by algebraic cycles (as seems reasonable given the results of \cite{kal} and the comments below), the answer is ``yes''. Furthermore, \cite[Section 5]{BO} conjectures that all crepant resolutions of $X$ have equivalent bounded derived categories of coherent sheaves. Confirmation of this conjecture would also give a positive answer.

In the particular case of $\xc$, it is possible to show that there
is a crepant resolution which can be described as a quiver
variety, \cite{nak}, generalising the $\mathbf{c} = \mathbf{0}$
case in \cite[Sections 1.3 and 1.4]{wang3} and the generic
$\mathbf{c}$ case in \cite[Section 11]{EG}. The $K$--theory of the
fibres has been studied, and is related to weight spaces of
integrable representations of Kac--Moody Lie algebras. Hence it is
reasonable to expect the number of simple
$\hc/\mathfrak{m}_x\hc$--modules also has this interesting
description. We will return to this in future work.


\begin{thebibliography}{10}

\bibitem{BO} A. Bondal and D. Orlov, Derived categories of coherent sheaves, {\it preprint}, AG/0206295.

\bibitem{BKR}
T. Bridgeland, A. King, and M. Reid,
The McKay correspondence as an equivalence of derived categories,
{\it J. Amer. Math. Soc.,} {\bf 14} (2001) 535--554.

\bibitem{BrGoo}
K.A.Brown and K.R.Goodearl, Homological aspects of noetherian PI Hopf algebras and irreducible modules of maximal dimension, {\it J.Algebra}, {\bf 198} (1997), 240--265.

\bibitem{BG}
K.A. Brown and I. Gordon, Poisson orders, representation theory and symplectic reflection algebras, {\it J. Reine Angew. Math.}, {\bf 559} (2003), 193--216.

\bibitem{CG}
N. Chriss and V. Ginzburg, Representation theory and complex geometry, Birkh\"{a}user, Boston, 1997.

\bibitem{CB}
W. Crawley-Boevey, Normality of Marsden-Weinstein reductions for
representations of quivers, {\it Math. Ann.,} {\bf 325} (2003),
55--79.

\bibitem{dCM}
M. de Cataldo and L. Migliorini,
The Chow motive of semismall resolutions, {\it preprint}, AG/0204067.

\bibitem{denloe}
J. Denef and F. Loeser, Germs of arcs on singular algebraic varieties and motivic integration, {\it Invent. Math.}, {\bf 135} (1999), 201--232.

\bibitem{EG}
P. Etingof and V. Ginzburg, Symplectic reflection algebras,
Calogero-Moser space, and deformed Harish-Chandra homomorphism,
{\it Invent. Math.,} {\bf 147} (2002) 243--348.

\bibitem{FG}
M. Finkelberg and V. Ginzburg, Calogero--Moser space and Kostka polynomials, {\it Adv.Math.}, {\bf 172} (2003), 137--150.

\bibitem{GK} V. Ginzburg and D. Kaledin, Poisson deformations of symplectic quotient singularities, {\it preprint},
AG/0212279, to appear {\it Adv.Math.}

\bibitem{gor} I. Gordon, Baby Verma modules for rational Cherednik algebras, {\it Bull.London Math.Soc.}, {\bf 35} (2003), 321--336.

\bibitem{hai} M. Haiman, Hilbert schemes, polygraphs and the
Macdonald positivity conjecture, {\it J. Amer. Math. Soc.}, {\bf
14} (4) (2001), 941--1006.

\bibitem{hai2} M. Haiman, Combinatorics, symmetric functions, and
Hilbert schemes, {\it preprint} (2003) to appear {\it Current
Developments in Mathematics, Proceedings of a conference at
Harvard University, November 2002}


\bibitem{hai3} M. Haiman, Cores, quivers and $n!$ conjectures, {\it in
preparation}

\bibitem{Ha}
R. Hartshorne, {\it Algebraic Geometry,} Graduate Texts in
Mathematics, no. 52, Springer-Verlag, 1977.

\bibitem{in} Y. Ito and I. Nakamura, McKay correspondence and Hilbert schemes, {\it Proc. Japan Acad. Ser. A Math. Sci.}, {\bf 72} (1996), 135--138.

\bibitem{kal0} D. Kaledin, Dynkin diagrams and crepant resolutions of quotient singularities, {\it preprint}, AG/9903157, to appear {\it Selecta Math.}

\bibitem{kal} D. Kaledin, McKay correspondence for symplectic quotient singularities, {\it Invent.Math.}, {\bf 148} (2002), no.1, 151--175.

\bibitem{KV}
M. Kapranov and E. Vasserot,
Kleinian singularities, derived categories and Hall algebras,
{\it Math. Ann.,} {\bf 316} (2000) 565--576.

\bibitem{K}
A.D. King,
Moduli of representations of finite--dimensional algebras,
{\it  Quart. J.  Math. Oxford Ser. (2),} {\bf 45} (1994) 515--530.

\bibitem{kraft}  H. Kraft, {\it Geometrische Methoden in der Invariantentheorie}, Aspects of Mathematics,
 D1. Friedr. Vieweg und Sohn, Braunschweig, 1984.

\bibitem{McCRob} J.C.McConnell and J.C.Robson, {\it Noncommutative Noetherian Rings}, Graduate studies in mathematics, 30 { American Mathematical Society}, 2001.

\bibitem{Mont}
S.M. Montgomery, {\em Fixed Rings of Finite Automorphism Groups of
Associative Rings,} Lecture Notes in Mathematics 818,
Springer-Verlag, 1980.


\bibitem{nak} H. Nakajima, Quiver varieties and Kac-Moody
algebras, \newblock {\it Duke Math. J.}, {\bf 91} (1998), no. 3, 515--560.

\bibitem{nakbook} H. Nakajima, {\it Lectures on Hilbert schemes of points on surfaces}, American Mathematical Society, Providence, RI, 1999.

\bibitem{vdb2} M. Van den Bergh, Non--commutative crepant
resolutions, {\it preprint}, RA/0211064, to appear {\it
Proceedings of the Abel Bicentennial Conference}.

\bibitem{ver} M.Verbitsky, Holomorphic symplectic geometry and orbifold singularities, {\it Asian J.Math.}, {\bf 4} (2000), no.3, 553--563.

\bibitem{wang2} W.Wang, Equivariant $K$-theory, wreath products, and Heisenberg algebra, {\it Duke Math. J.}, {\bf 103} (2000), no. 1, 1--23.

\bibitem{wang3} W.Wang, Hilbert schemes, wreath products, and the McKay
correspondence, \textit{preprint}, math.AG/9912104

\bibitem{wang1} W.Wang, Algebraic structures behind Hilbert schemes and wreath products, in \textit{Recent developments in infinite-dimensional Lie algebras and conformal field theory} (Charlottesville, VA, 2000), 271--295, Contemp. Math., 297, Amer. Math. Soc., 2002.

\end{thebibliography}
\end{document}